\pdfoutput=1
\RequirePackage{ifpdf}
\ifpdf 
\documentclass[pdftex]{sigma}
\else
\documentclass{sigma}
\fi

\usepackage{amscd,array,stmaryrd}
\usepackage[all]{xy}

\numberwithin{equation}{section}

\newtheorem{Theorem}{Theorem}[section]

\newtheorem{Lemma}[Theorem]{Lemma}
\newtheorem{Proposition}[Theorem]{Proposition}
\newtheorem{Conjecture}[Theorem]{Conjecture}
{\theoremstyle{definition}
\newtheorem{Definition}[Theorem]{Definition}

\newtheorem{Remark}[Theorem]{Remark}
}

\def\a{\alpha}

\def\g{\gamma}
\def\Ga{\Gamma}
\def\d{\delta}

\def\l{\lambda}

\def\m{\mu}

\newcommand{\bbD}{\mathbb{D}}
\newcommand{\bbG}{\mathbb{G}}

\newcommand{\bbR}{\mathbb{R}}

\newcommand{\bbN}{\mathbb{N}}

\newcommand{\bbQ}{\mathbb{Q}}

\newcommand{\aff}{\mathfrak{aff}}

\newcommand{\cC}{{\mathcal{C}}}

\newcommand{\ce}{\mathfrak{ce}}
\newcommand{\cf}{\mathfrak{cf}}

\newcommand{\cU}{{\mathcal{U}}}

\newcommand{\D}{\mathcal{D}}

\newcommand{\Dlm}{\mathcal{D^{\l,\m}}}

\newcommand{\cE}{\mathcal{E}}

\newcommand{\Div}{\mathrm{Div}}

\newcommand{\cF}{{\mathcal{F}}}

\newcommand{\fkg}{{\mathfrak{g}}}
\newcommand{\gl}{{\mathrm{gl}}}

\newcommand{\Hom}{\mathrm{Hom}}

\newcommand{\fkh}{{\mathfrak{h}}}
\newcommand{\rH}{\mathrm{H}}

\newcommand{\Id}{\mathrm{Id}}
\newcommand{\im}{\mathrm{im}}
\newcommand{\ccL}{\mathcal{L}}

\newcommand{\LD}{\mathcal{L}^{\lambda,\mu}}
\newcommand{\LDX}{\mathcal{L}_X^{\lambda,\mu}}
\newcommand{\LDi}{\mathcal{L}_{X_i}^{\lambda,\mu}}
\newcommand{\Ll}{\ell_X^\lambda}
\newcommand{\Lm}{\ell_X^\mu}

\newcommand{\Ld}{L^\delta}
\newcommand{\bLd}{\bar{L}^{-\delta}}
\newcommand{\LdX}{L_X^\delta}

\newcommand{\Ldi}{L_{X_i}^\delta}

\newcommand{\ro}{\mathfrak{o}}

\newcommand{\pol}{\mathrm{pol}}

\newcommand{\cQ}{{\mathcal{Q}}}
\newcommand{\cQlm}{\mathcal{Q}^{\l,\m}}

\newcommand{\Qlm}{Q^{\l,\m}}

\newcommand{\cN}{{\mathcal{N}}}

\newcommand{\cS}{{\mathcal{S}}}
\newcommand{\bcS}{{\bar{\mathcal{S}}}}

\newcommand{\SL}{\mathrm{SL}}

\newcommand{\Sl}{\mathfrak{sl}}

\newcommand{\Vect}{\mathrm{Vect}}

\newcommand{\half}{\frac{1}{2}}

\begin{document}

\allowdisplaybreaks

\renewcommand{\PaperNumber}{022}

\FirstPageHeading

\ShortArticleName{Conformally Equivariant Quantization~-- a~Complete Classif\/ication}

\ArticleName{Conformally Equivariant Quantization~--\\ a~Complete Classif\/ication}

\Author{Jean-Philippe MICHEL}
\AuthorNameForHeading{J.-P.~Michel}
\Address{University of Luxembourg, Campus Kirchberg, Mathematics Research Unit,\\
 6, rue Richard Coudenhove-Kalergi, L-1359 Luxembourg City, Luxembourg}
\Email{\href{mailto:jean-philippe.michel@uni.lu}{jean-philippe.michel@uni.lu}}
\URLaddress{\url{http://math.uni.lu/~michel/SitePerso/indexLux.html}}

\ArticleDates{Received July 29, 2011, in f\/inal form April 11, 2012; Published online April 15, 2012}

\Abstract{Conformally equivariant quantization is a peculiar map between symbols of real weight~$\delta$ and dif\/ferential operators acting on tensor densities, whose real weights are designed by $\lambda$ and $\lambda+\delta$. The existence and uniqueness of such a map has been proved by Duval, Lecomte and Ovsienko for a generic weight $\delta$.
Later, Silhan has determined the critical values of $\delta$ for which unique existence is lost, and conjectured that for those values of $\delta$ existence is lost for a generic weight $\lambda$. We fully determine the cases of existence and uniqueness of the conformally equivariant quantization in terms of the values of $\delta$ and $\lambda$. Namely, (i)~unique existence is lost if and only if there is a nontrivial conformally invariant dif\/ferential operator on the space of symbols of weight $\delta$, and (ii)~in that case the conformally equivariant quantization exists only for a f\/inite number of $\lambda$, corresponding to nontrivial conformally invariant dif\/ferential operators on $\lambda$-densities. The assertion~(i) is proved in the more general context of IFFT (or AHS) equivariant quantization.}

\Classification{53A55; 53A30; 17B56; 47E05}

\section{Introduction}

Quantization originates from physics, and lies at the correspondence between classical and quantum formalisms. It has given the impetus for the development of numerous mathematical theo\-ries, and admits as many dif\/ferent def\/initions.
In this paper, by quantization we mean the inverse of a certain symbol map, i.e.\ a linear map from functions on a cotangent bundle $T^*M$, which are polynomial in the f\/iber variables, to the space of dif\/ferential operators on $M$. No such quantization can be canonically def\/ined from the dif\/ferential geometry of $M$. The idea of equivariant quantization is to build one from a richer geometric structure on $M$: the local action of a Lie group $G$, or equivalently, of its Lie algebra $\fkg$. This imposes the condition that the manifold $M$ possesses a locally f\/lat  $G$-structure.

Equivariant quantization was f\/irst developed for projective and conformal structures by Duval, Lecomte and Ovsienko \cite{DLO99,LOv99}. There, the condition of equivariance singles out a unique quantization up to normalization. Generalizations to the curved case have been proposed recently in terms of Cartan connections \cite{MRa07,MRa10b} or tractors \cite{CSi09,Sil09}, but we restrict ourselves here to the locally f\/lat case, or equivalently to $\bbR^n$.
In \cite{BMa06}, Boniver and Mathonet exhibit the correct class of Lie algebras to be considered for equivariant quantization. It is provided by maximal Lie subalgebras $\fkg$ among those of f\/inite dimension in $\Vect_\pol(\bbR^n)$, the Lie algebra of polynomial vector f\/ields. They prove in \cite{BMa01} that these Lie algebras are precisely the \textit{irreducible filtered Lie algebras of finite type} (with no complex structure) classif\/ied by Kobayashi and Nagano~\cite{KNa64}, referred to as IFFT-algebras. Particular examples are the Lie algebra ${\Sl_{n+1}\simeq\Sl(n+1,\bbR)}$ of projective vector f\/ields, and the Lie algebra $\cf\simeq\ro(p+1,q+1)$ of conformal Killing vector f\/ields on $(\bbR^n,\eta)$, where $n=p+q$ and $\eta$ is the f\/lat metric of signature $(p,q)$. Generali\-zing~\cite{BMa06}, Cap and Silhan prove in~\cite{CSi09} the existence and uniqueness of the $\fkg$-equivariant quantization with values in dif\/ferential operators acting on irreducible homogeneous bundles, barring certain exceptional bundles. The determination of these exceptional cases and their geometric interpretation is an open problem, since they appear in the seminal work~\cite{DLO99} dealing with bundles of tensor densities $|\Lambda^n T^*\bbR^n|^{\otimes\l}$ of arbitrary weight $\l\in\bbR$. Our paper is devoted to a step towards its solution. We give a complete resolution in the original case of conformally equivariant quantization ${\cQlm:\cS^\d\rightarrow\Dlm}$ with values in the space $\D^{\l,\m}$ of dif\/ferential operators from $\l$- to $\m$-tensor densities, the symbols being here of weight $\d=\m-\l$.

\looseness=1
This work can be regarded in the broader perspective of invariant bidif\/ferential operators, whose equivariant quantization is an example. Another example is provided by generalized transvectants (or Rankin--Cohen bracket), whose arguments are spaces of weighted densities. There the same phenomenon occurs: they exist and are unique except for certain exceptional weights \cite{ORe03}. In fact, this is general, as proved by Kroeske in his thesis~\cite{Kro08a}, where he stu\-dies invariant bidif\/ferential operators for parabolic geometries. From a number of examp\-les and a~proof in the projective case~\cite{Kro08}, Kroeske proposes a paradigm which can be roughly phrased as follows: every exceptional case of existence or uniqueness of an invariant bidif\/ferential  ope\-ra\-tor originates from the existence of a nontrivial invariant dif\/ferential operator on one of the factors. We show that this is particularly illuminating in the case of equivariant quantization. Namely, resorting to the interpretation in terms of cohomology of $\fkg$-modules deve\-lo\-ped by Lecomte in~\cite{Lec00}, we prove our f\/irst main result: unique existence of $\fkg$-equivariant quantization is lost if and only if there exists a $\fkg$-invariant operator on the space of symbols. Remarkably, this correspondence is proved directly rather than obtained \`a posteriori, like in the work of Kroeske. The theory of invariant operators being well-developed, we get an ef\/f\/icient way to determine the exceptional bundles for equivariant quantization. Retur\-ning to the aforementioned conformally equivariant quantization ${\cQlm:\cS^\d\rightarrow\Dlm}$, we recover the set of {\it critical} values of $\d$ for which unique existence is lost, already determined by Silhan~\cite{Sil09}.

\looseness=1
The next natural question is to determine for which \textit{resonant} pairs of irreducible homogeneous bundles the $\fkg$-equivariant quantization exists but is not unique. This turns to be an harder question, and we address it only for pairs of line bundles of tensor densities. By our f\/irst main theorem, we can obtain the critical shift $\d=\m-\l$ of their weights for which there is not existence and uniqueness of the $\fkg$-equivariant quantization. We have then to determine the resonant values $(\l,\m)$ for which the $\fkg$-equivariant quantization still exists. In the projective case, an explicit formula for the quantization provides the answer \cite{DOv01a, LOv99}. Moreover, the situation has been fully understood in terms of cohomology of $\Sl_{n+1}$-modules, the existence of the projectively equivariant quantization being characterized by the triviality of a $1$-cocycle, which depends on~$\d$ and~$\l$~\cite{Lec00}. Using a similar approach, we obtain our second main result: for each critical value~$\d$ of the conformally equivariant quantization $\cQlm$, there is a f\/inite number of resonances $(\l,\l+\d)$ that we determine. This completes the known results on symbols of degree at most~$3$, and proves a conjecture of Silhan: the conformally equivariant quantization does generically not exist for~$\d$ critical. In addition, we provide an interpretation along the lines of the Kroeske's paradigm: the resonances correspond to conformally invariant operators on the source space of densities, e.g.\ the conformal powers of the Laplacian. This allows us to construct conformally equivariant quantization in the resonant cases. Together with Silhan's work \cite{Sil09}, this provides construction for the conformally equivariant quantization whenever it exists.

\looseness=1
Let us outline the contents of this paper. Section~\ref{section2} is devoted to the proof of Theorem~\ref{Thmgequi}: for $\fkg$ an IFFT-algebra, the $\fkg$-equivariant quantization exists and is unique if and only if there is no $\fkg$-invariant dif\/ferential operator on the space of symbols which strictly lo\-wers the degree. It relies on the fact that equivariant quantization exists if a certain $1$-cocycle is a coboun\-da\-ry. We know from \cite{BMa06} that this is the case for generic $\d'$, and this property is stable in the limit $\d'\rightarrow\d$ if the space of $1$-coboundaries is of constant dimension. That condition happens to be equivalent to the absence of $\fkg$-invariant operators on $\d$-weighted symbols. In Sections~\ref{section3} and~\ref{section4} we illustrate and complete Theorem~\ref{Thmgequi} for, respectively, $\fkg$ the projective and the conformal Lie algebras, the quantization being valued in dif\/ferential ope\-ra\-tors acting on tensor densities. First, we determine the critical values via the classif\/ication of invariant operators on the space of symbols. Then we prove that those invariant operators give rise to nontrivial $1$-cocycles. They obstruct existence of the equivariant quantization except when there is an invariant dif\/ferential operator from $\l$-densities to a certain homogeneous bundle.  Finally, the latter operator allows us to construct  an equivariant quantization, proving thus its existence for exactly those values of~$\l$. That leads to the complete list of resonances $(\l,\l+\d)$, together with their interpretation in terms of invariant operators on the space of $\l$-densities. We end Section~\ref{section4} with a detailed treatment of the symbols of degree less than~$3$,  thus interpreting the results of Loubon~Djounga~\cite{Lou01} along the line of the Kroeske's paradigm. The last section gives us the opportunity to propose some natural extensions of our results.

Throughout this paper, the space of linear maps between two real vector spaces $V$ and $W$ is denoted by $\Hom(V,W)$, and its elements are called operators. We work on~$\bbR^n$ and the Einstein summation convention is understood.

\section{On the existence and uniqueness of equivariant quantization}\label{section2}

\subsection{Def\/inition of equivariant quantization}\label{section2.1}

We start with the def\/initions of the algebra $\D(\bbR^n)$ of dif\/ferential operators on $\bbR^n$, and its algebra of symbols $\cS(\bbR^n)$. The former is f\/iltered by the subspaces $\D_k(\bbR^n)$ of dif\/ferential operators of order $k$, def\/ined as the spaces of operators $A$ on $\cC^\infty(\bbR^n)$ satisfying $[\dots [A,f_0],\dots],f_k]=0$ for all functions $f_0,\ldots,f_k\in\cC^\infty(\bbR^n)$, considered here as (zero order) operators on $\cC^\infty(\bbR^n)$. The latter is the canonically associated graded algebra, def\/ined by
$\cS(\bbR^n)=\bigoplus_{k=0}^\infty\D_k(\bbR^n)/\D_{k-1}(\bbR^n)$. This may be identif\/ied with the algebra of functions on $T^*\bbR^n$ that are polynomial in the f\/ibers, the grading corresponding to the polynomial degree. The canonical projection $\sigma_k:\D_k(\bbR^n)\rightarrow\D_k(\bbR^n)/\D_{k-1}(\bbR^n)$ is called the principal symbol map. It admits a section, the normal ordering, given by
\begin{gather}\label{normalordering}
\cN: \ P^{i_1\dots i_k}(x)p_{i_1}\cdots p_{i_k}\mapsto P^{i_1\dots i_k}(x)\partial_{i_1}\cdots \partial_{i_k},
\end{gather}
where $(x^i,p_i)$ are coordinates on $T^*\bbR^n$. This def\/ines a linear isomorphism $\cS(\bbR^n)\simeq\D(\bbR^n)$.

We are interested in the action of vector f\/ields on these both algebras. First of all, we intro\-duce the $\Vect(\bbR^n)$-module of $\l$-densities $\cF^\l=(\cC^\infty(\bbR^n),\ell^\l)$ as a one parameter deformation of the module $\cC^\infty(\bbR^n)$, the $\Vect(\bbR^n)$-action being given by $X\mapsto\ell_X^\l=X+\l\Div(X)$, with $\Div$ the divergence and $\l\in\bbR$ the weight of the densities. This module corresponds geometrically to the space of sections of the trivial line bundle $|\Lambda^n T^*\bbR^n|^{\otimes\l}$. It gives rise to the $\Vect(\bbR^n)$-module $\Dlm=(\D(\bbR^n),\LD)$ of dif\/ferential operators from $\l$- to $\mu$-densities, endowed with the adjoint action
\[\LDX A=\Lm A-A\Ll,\]
for all $X\in\Vect(\bbR^n)$ and $A\in\D(\bbR^n)$. This action preserves the f\/iltration of $\D(\bbR^n)$, and so the algebra of symbols inherits a $\Vect(\bbR^n)$-module structure compatible with the grading. We denote this structure by $\cS^\d=(\cS(\bbR^n),\Ld)$, where $\d=\mu-\l$ is the shift. The action of $\Vect(\bbR^n)$ is given in coordinates by
\begin{gather}\label{LdX}
\LdX=X^i\partial_i-p_j(\partial_iX^j)\partial_{p_i}+\d\,\Div(X),
\end{gather}
and coincides with the canonical action on functions on $T^*M$ tensored with $\d$-densities. We denote by $\cS^\d_k$ the submodule of homogeneous symbols of degree $k$. 
\begin{Definition}
Let $\fkg$ be a Lie subalgebra of $\Vect(\bbR^n)$. A $\fkg$-equivariant quantization is a~$\fkg$-module morphism
\[
\cQlm: \ \cS^\d\rightarrow\Dlm,
\]
such that $\cQlm$ is a right inverse of the principal symbol map on homogeneous symbols.
\end{Definition}
Using the normal ordering \eqref{normalordering}, we freely identify $\Dlm$ with $(\cS(\bbR^n),\LD)$, where we keep the same notation for $\LD$ and its pull-back on symbols by $\cN$.
It is then a matter of computation to prove that $\LDX=\LdX$ if $X$ is an af\/f\/ine vector f\/ield and $\d=\m-\l$. In other words, $\cN$ is an $\aff(n,\bbR)$-equivariant quantization~\cite{DLO99}.

\subsection{IFFT-algebras and equivariant quantization}\label{section2.2}

\looseness=1
An IFFT-algebra $\fkg$ is a simple Lie subalgebra of the polynomial vector f\/ields on $\bbR^n$. As such, it admits a gradation  by the degree of the vector f\/ield components  $\fkg=\fkg_{-1}\oplus\fkg_0\oplus\fkg_1$ which is compatible with the bracket: $[\fkg_i,\fkg_j]=\fkg_{i+j}$, where $\fkg_i=\{0\}$ if $i\notin\{-1,0,1\}$. The Lie subalgebra $\fkg_{-1}$ consists of translations, $\fkg_1$ consists of so-called $\fkg$-inversions, and $\fkg_0$ contains the dilation and acts irreducibly on both $\fkg_{-1}$ and $\fkg_1$. Consequently, $\fkg$-invariance is equivalent to $\fkg_{-1}\oplus\fkg_0$-invariance plus the invariance with respect to one non-zero element in~$\fkg_{1}$.

Let us introduce some notation. For $\fkh$ a Lie subalgebra of $\Vect(\bbR^n)$, the space of $\fkh$-invariant operators between two $\fkh$-submodules $F\subset\cS^\d$ and $F'\subset\cS^{\d'}$ is def\/ined as
\begin{gather}\label{Homfkh}
\Hom_\fkh(F,F')=\{A\in\Hom(F,F')\,|\;\forall\, X\in\fkh,\, [L_X^*,A]=0\},
\end{gather}
where the commutator $[L_X^*,A]$ is a symbolic notation for $L_X^{\d'}A-A\LdX$. Note that the vector f\/ields which preserve $F$ and $F'$ clearly act on $\Hom_\fkh(F,F')$. We def\/ine similarly the space of $\fkh$-invariant dif\/ferential operators $\D_\fkh(F,F')$ and its module structure.
\begin{Lemma}[\cite{LOv99}]\label{E}
Let $\fkg$ be an IFFT-algebra, and fix $l< k$. The vector space $\Hom_{\fkg_{-1}\oplus\fkg_0}(\cS^{\d}_k,\cS^{\d}_l)$ is finite-dimensional, independent of $\d$, and equal to the space $\D_{\fkg_{-1}\oplus\fkg_0}(\cS^{\d}_k,\cS^{\d}_l)$ of invariant differential operators.
\end{Lemma}
\begin{proof}
We only need the fact that $\fkg_{-1}\oplus\fkg_0$ contains the translations and the dilation. Indeed, an operator $A\in\Hom(\cS^\d_k,\cS^\d_l)$ invariant with respect to those transformations was proved in~\cite{LOv99} to be a dif\/ferential operator if $l< k$. Hence it is generated by the coordinates $p_i$ and the derivatives~$\partial_i$,~$\partial_{p_i}$. Now the dilation invariance ensures that the degree in $\partial_i$ is equal to the degree in $\partial_{p_i}$ minus the degree in~$p_i$. The degree in $\partial_{p_i}$ being $k$ at most, the space $\Hom_{\fkg_{-1}\oplus\fkg_0}(\cS^\d_k,\cS^\d_l)$ is f\/inite-dimensional. Since the divergence of af\/f\/ine vector f\/ields is a constant, the action of $\fkg_{-1}\oplus\fkg_0$ on $\Hom(\cS^{\d}_k,\cS^{\d}_l)$ is independent of $\d$ and hence the subspace of $\fkg_{-1}\oplus\fkg_0$-invariant operators is also.
\end{proof}

We turn now to the study of the $\fkg$-equivariance condition for a quantization~$\cQlm$. We can always factor the latter through the normal ordering, def\/ining the linear automorphism~$\Qlm$ of~$\cS^\d$ by~$\cQlm=\cN\circ\Qlm$. Restricting to~$\cS^\d_k$, we then get
\begin{gather}\label{Qlmphi}
\Qlm=\sum_{l=0}^k\phi_l,
\end{gather}
where $\phi_0=\Id$ and $\phi_l\in\Hom(\cS^\d_k,\cS^\d_{k-l})$. Since the normal ordering is $\aff(n,\bbR)$-equivariant, we deduce from Lemma~\ref{E} that more precisely $\phi_l\in\D_{\fkg_{-1}\oplus\fkg_0}(\cS^\d_k,\cS^\d_{k-l})$ for all $l=1,\ldots,k$. The full equivariance condition for $\cQlm$ reads then on $\cS^\d_k$ as
\begin{gather}\label{Qequiv}
\big[\phi_l,\LdX\big]=\big(\LDX-\LdX\big)\phi_{l-1},
\end{gather}
for all $l=1,\ldots,k$, where $X$ is a non-zero element of $\fkg_1$. To solve it directly is too intricate except for $\fkg=\Sl_{n+1}$. Nevertheless, the following theorem has been proven in \cite{BMa06}, resorting to simultaneous diagonalization of the Casimir operators on the modules of symbols and dif\/ferential operators. We also refer to \cite{CSi09}, where the additional hypothesis of simplicity of the semi-simple part of $\fkg_0$ is dropped.

\begin{Theorem}[\cite{BMa06,CSi09}]\label{Genericthm}
Let $\fkg$ be an IFFT-algebra and fix $k\in\bbN$. The $\fkg$-equivariant quantization $\cQlm:\cS^\d_k\rightarrow\D^{\l,\m}_k$ exists and is unique if $\m-\l=\d\notin I_k$, with $I_k$ a finite subset of $\bbQ$.
\end{Theorem}

\subsection{Equivariant quantization and Lie algebra cohomology}\label{section2.3}
We give here a cohomological interpretation of the equations \eqref{Qequiv} encoding the equivariance of $\Qlm$.
Let us f\/irst give a brief review of the cohomology of $\fkg$-modules and its link with the splitting of exact sequences of $\fkg$-modules (see \cite{Fuk87} for more details). The cohomology of the $\fkg$-module $(M,L^M)$ is def\/ined in terms of the $k$-cochains, which are the linear maps $\Lambda^k \fkg\rightarrow M$, and the dif\/ferential $d$ which reads, on $0$- and $1$-cochains $\phi$ and $\g$,
\begin{gather}
d\phi(X) =L^M_X\phi,\qquad \label{1cocycle}
d\g(X,Y) =L^M_X\g(Y)-L^M_Y\g(X)-\g([X,Y]),
\end{gather}
where $X$, $Y$ are in $\fkg$. Let now $(A,L^A)$, $(B,L^B)$, $(C,L^C)$ be three $\fkg$-modules, and suppose that we have an exact sequence of $\fkg$-modules,
\[
\xymatrix{
0\ar[r]& (A,L^A)\ar[r]^{\iota} &(B,L^B) \ar[r]^\sigma &(C,L^C)\ar @/^1pc/@{.>}[l]^{\tau} \ar[r] & 0,
}
\]
with $\tau$ a linear section. This def\/ines a $1$-cocycle $\g=\iota^{-1}(L^B\circ\tau-\tau\circ L^C)$ with values in $\Hom(C,A)$. Its cohomolgical class does not depend on the choice of linear section $\tau$, and the sequence of $\fkg$-modules is split if and only if $\g=d\phi$ is a coboundary. The splitting morphism is then $\tau+\iota\circ\phi$. Moreover, if $\g$ vanishes on a Lie subalgebra $\fkh$ then $\phi$ is $\fkh$-invariant.

The existence and uniqueness of a $\fkg$-equivariant quantization can be rephrased in terms of cohomology of $\fkg$-modules \cite{Lec00}. Indeed, such a quantization exists if for every $k\in\bbN$ the exact sequence
\[
0\rightarrow\D^{\l,\m}_{k-1}\rightarrow\D^{\l,\m}_{k}\rightarrow \cS^\d_k \rightarrow 0
\]
is split in the category of $\fkg$-modules. Using the normal ordering as a linear splitting, this means that the $1$-cocycle $\g=\LD-\Ld$ admits a trivial cohomology class $[\g]$. By vanishing of $\g$ on the af\/f\/ine part of $\fkg$, the latter pertains to the following relative cohomology space $H^1(\fkg,\fkg_{-1}\oplus\fkg_0;\Hom(\cS^\d_k,\D^{\l,\m}_{k-1}))$.
 The modules $\Hom(\cS^\d_k,\D^{\l,\m}_{k-1})$ are quite complex to handle, but modding out by $\D^{\l,\m}_{k-l}$ for increasing $l$, we are reduced by induction to the simpler modules $\Hom(\cS^\d_k,\cS^\d_{k-l})$. Thus, a $\fkg$-equivariant quantization on $F^\d$, a submodule of $\cS^\d_k$, is a section of $\fkg$-modules $\psi_k:F^\d\rightarrow\D^{\l,\m}_k$ def\/ined inductively by $\psi_0=\Id$ and the commutative triangle in the following diagram of $\fkg$-modules
\begin{gather}\label{exactseq}
\begin{split}
\xymatrix{
0\ar[r] & \cS^\d_{k-l}\ar[r] & \D^{\l,\m}_{k}/\D^{\l,\m}_{k-l-1}\ar[r] & \D^{\l,\m}_{k}/\D^{\l,\m}_{k-l} \ar[r] & 0\\
&&& F^\d \ar@{^{(}->}[u]_{\psi_{l-1}}\ar@{^{(}->}[ul]^{\exists ? \psi_{l}} &
}
\end{split}
\end{gather}
for successively all $l=1,\ldots,k$. With notation as in \eqref{Qlmphi}, we have $\psi_l=\sum_{i=0}^l\phi_i$. In the next lemma we recover the equivariance condition \eqref{Qequiv} using this cohomological approach.
\begin{Lemma}\label{lemcoh}
The partial quantization $\psi_l$ defined in \eqref{exactseq} exists if and only if $\psi_{l-1}$ exists and the $1$-cocycle $\g_l=(\LD-\Ld)\phi_{l-1}$, whose class belongs to $H^1(\fkg,\fkg_{-1}\oplus\fkg_0;\Hom(F^\d,\cS^\d_{k-l}))$, is trivial. In this case, we have $\psi_l=\psi_{l-1}+\phi_l$, where $\phi_l$ satisfies
\begin{gather}\label{phig}
\big[\phi_{l},\LdX\big]=\g_l(X)
\end{gather}
for some non-zero $X$ in $\fkg_1$ and $\phi_{l}\in\D_{\fkg_{-1}\oplus\fkg_0}(F^\d,\cS^\d_{k-l})$.
\end{Lemma}
\begin{proof}
Via the normal ordering, $\psi_{l-1}$ lifts linearly to $\D^{\l,\m}_{k}/\D^{\l,\m}_{k-l-1}$. The existence of the morphism $\psi_{l}$ relies then on the triviality of the $1$-cocycle $(\LD-\Ld)\psi_{l-1}$. As it takes values in $\cS^\d_{k-l}$, and as $\LDX-\LdX$ lowers the degree by one, the previous $1$-cocycle is equal to $\g_l$.
Thus $\psi_l$ exists if and only if $\g_l=d\phi_{l}$, and $\psi_l=\psi_{l-1}+\phi_l$ is the splitting morphism. From the def\/inition of $\g_l$, we get the $\fkg_{-1}\oplus\fkg_0$-invariance of $\phi_l$. Finally, the irreducibility of the action of $\fkg_0$ on $\fkg_1$ together with def\/inition of $1$-cocycles shows that \eqref{phig} is equivalent to the triviality of $\g_l$.
\end{proof}

\subsection{Main result}\label{section2.4}

 We can now give a characterization of the critical values $\d$ of the $\fkg$-equivariant quantization in terms of $\fkg$-invariant operators.
\begin{Theorem}\label{Thmgequi}
Let $\fkg$ be an IFFT-algebra and let $F^\d=(F,\Ld)$ be a $\fkg$-submodule of $\cS^\d_k$ for every~\mbox{$\d\in\bbR$}. The $\fkg$-equivariant quantization exists and is unique on $F^\d$ if and only if there exists no non-zero $\fkg$-invariant differential operator from $F^\d$ to $\cS^\d_{k-l}$, for $l=1,\ldots,k$.
\end{Theorem}
\begin{proof}
Let $\d\in I_k$, def\/ined in Theorem \ref{Genericthm}. Clearly, if there exists a $\fkg$-equivariant quantization on $F^\d$, its uniqueness is equivalent to the absence of $\fkg$-invariant operator from $F^\d$ to $\cS^\d_{k-l}$ for $l=1,\ldots,k$. We only have to prove that such an absence implies existence of the $\fkg$-equivariant quantization on $F^\d$. By induction on $l$, this amounts to obtaining the partial quantization $\psi_l$ out of $\psi_{l-1}$. By the preceding lemma, this means to prove that the $1$-cocycle $\g_l^\d$ is trivial: $\g^\d_l(X)\in[E_0^\d,L^\d_X]$ for a non-zero $X$ in $\fkg_1$ and $E_0^\d:=\D_{\fkg_{-1}\oplus\fkg_0}(\cS^{\d}_k,\cS^{\d}_{k-l})$.

By Theorem~\ref{Genericthm}, we know that $\fkg$-equivariant quantization exists for shifts $\d'\neq\d$ in a small enough neighborhood $\cU$ of $\d$. In particular this implies that $\g_l^{\d'}(X)$ is a coboundary, hence it pertains to the space $[E_0^{\d'},L_{X}^{\d'}]$. We have to show that this remains true in the limit $\d'\rightarrow\d$. The Lemma \ref{E} ensures that the domains $E_0^{\d'}$ are f\/inite-dimensional and independent of~$\d'$, so we denote them all by $E_0$. We also introduce $E_1$, the subspace of $\D((\cS(\bbR^n),\cS(\bbR^n))$ generated by the family of spaces $ [E_0,L_X^{\d'} ]$ for $\d'\in\bbR$. As $X$ is quadratic, we deduce from \eqref{LdX} that the space generated by the operators $L_X^{\d'}$ is f\/inite-dimensional, and hence $E_1$ is also. Consequently, we get a continuous family of linear maps, indexed by $\d'\in\bbR$, between f\/inite-dimensional spaces:
\[
 \big[\cdot,L_{X}^{\d'}\big]: \ E_0
 \longrightarrow E_1.
\]
 Since  there is no $\fkg$-invariant dif\/ferential operator from $F^{\d'}$ to $\cS^{\d'}_{k-l}$ for $\d'\in\cU$, the kernel of $[\cdot,L_{X}^{\d'}]$ is reduced to zero  for $\d'\in\cU$. Consequently, the spaces $\im([\cdot,L_{X}^{\d'}])$ are of constant rank on $\cU$ and
the relation $\g^{\d'}_l(X)\in\im([\cdot,L_{X}^{\d'}])$ is preserved in the limit $\d'\rightarrow \d$.
\end{proof}
The latter proof does not give a direct construction of the $\fkg$-equivariant quantization, but it completes its usual construction in terms of Casimir operators \cite{BMa06}. Indeed, if the latter method fails for a shift $\d$ whereas the $\fkg$-equivariant quantization $\cQ^{\l,\l+\d}$ exists and is unique, we have just shown that it is recovered from $\cQ^{\l,\l+\d'}$ in the limit $\d'\rightarrow\d$.

Let $V$ and $W$ be irreducible representations of $\fkg_0$ of f\/inite dimensions, and let $\mathcal{V}=\bbR^n\times V$, $\mathcal{W}=\bbR^n\times W$ the corresponding trivial bundles over $\bbR^n$. A natural generalization of $\fkg$-equivariant quantization is to consider dif\/ferential operators from $\Ga(\mathcal{V})$ to $\Ga(\mathcal{W})\otimes\cF^\d$, the space of symbols being then $\cS^\d\otimes\Ga(\mathcal{V}^*)\otimes\Ga(\mathcal{W}):=\cS^\d(\mathcal{V},\mathcal{W})$.
Theorem \ref{Genericthm} has been generalized in \cite{CSi09} to this context, and all the results of this section generalize straightforwardly to that situation also (for Lemma \ref{E}, note that the dilation acts diagonally on $V$ and $W$). This leads to the following theorem.
\begin{Theorem}\label{ThmgequiAHS}
Let $\fkg$ be an IFFT-algebra and let $F^\d=(F,\Ld)$ be a $\fkg$-submodule of $\cS^\d_k(\mathcal{V},\mathcal{W})$ for any $\d\in\bbR$. The $\fkg$-equivariant quantization exists and is unique on $F^\d$ if and only if there exists no non-zero $\fkg$-invariant differential operator from $F^\d$ to $\cS^\d_{k-l}(\mathcal{V},\mathcal{W})$, for $l=1,\ldots,k$.
\end{Theorem}
Let us mention that we obtain a necessary and suf\/f\/icient condition for $\d$ to be a critical value, contrary to the previous works relying on the diagonalization of the Casimir operator on the space of symbols. The suf\/f\/icient condition obtained there was that specif\/ied eigenvalues of this operator are equal. This is clearly a stronger condition on $\d$ than ours.

\section{Projectively equivariant quantization}\label{section3}
We turn now to the case $\fkg=\Sl_{n+1}$ and restrict our consideration to dif\/ferential operators acting on densities. First, we recall the construction of an explicit formula for the projectively equivariant quantization, using freely results of the original works \cite{DOv01a, LOv99}. Then we study in detail the critical and resonant values, in particular their link with existence of projectively invariant dif\/ferential operators. This can be seen as a warm up for the conformal case.

\subsection{Explicit formula}\label{section3.1}
The projective action of $ \SL(n+1,\bbR)$ on the projective space $\bbR P^n$ induces an embedding of $\Sl(n+1,\bbR)$ into the polynomial vector f\/ields on $\bbR^n$. The resulting Lie algebra $\Sl_{n+1}$ of projective vector f\/ields is generated by the af\/f\/ine vector f\/ields and the projective inversions
\[
X^i=x^ix^j\partial_j,
\]
for $i=1,\ldots,n$. This shows that $\Sl_{n+1}$ is an IFFT-algebra. Recall that for af\/f\/ine vector f\/ields $X$, we have $\LDX=\LdX$, so the lack of projective equivariance for the normal ordering is described on $\cS^\d_k$ by
\begin{gather}\label{ProjLDLd}
\ccL^{\l,\m}_{X^i}-L^\d_{X^i}=\ell_{k-1}(\l)\partial_{p_i},
\end{gather}
where $\ell_k(\l)=-(k+\l (n+1))$. Resorting to the preceding section, the projectively equivariant quantization on $\cS^\d_k$ decomposes as $\cN\circ\left(\Id+\phi_{1}+\dots+\phi_{k}\right)$, where $\phi_m$ is an $\aff(n,\bbR)$-invariant operator lowering the degree by $m$ and satisfying
\begin{gather}\label{CocycleProj}
\big[\phi_{m},L^\d_{X^i}\big]=\g_m\big(X^i\big),
\end{gather}
with $\g_m=(\LD-\Ld)\phi_{m-1}$. Weyl's theory of invariants \cite{Wey97}, together with Lemma~\ref{E}, shows that the $\aff(n,\bbR)$-invariant dif\/ferential operators acting on symbols are generated by the Euler operator $\cE=p_i\partial_{p_i}$ and the divergence operator $D=\partial_i\partial_{p_i}$. Hence, restricted to~$\cS^\d_k$, the map~$\phi_{m}$ is of the form~$c^k_mD^{m}$, with $c^k_m\in\bbR$. This coef\/f\/icient is determined by substitution into the equation \eqref{CocycleProj} and using
\begin{gather}\label{Dproj}
\big[D^m,L^\d_{X^i}\big]|_{\cS^\d_k}=d^k_m(\d)\partial_{p_i}D^{m-1},
\end{gather}
where $d^k_m(\d)=m(-2k+m+1+(\d-1) (n+1))$.
\begin{Theorem}[\cite{DOv01a, LOv99}]\label{QEP}
Let $1\leq l\leq k$. If $\d\neq 1+\frac{2k-l-1}{n+1}$, there exists a unique projectively equivariant quantization on $\cS^\d_k$, given by $\cN\circ\left(\sum_{m=0}^{k}c^k_m D^m\right)$, with $c^k_0=1$ and
\[
c^k_{m}=\frac{\ell_{k-m}(\l)}{d^k_m(\d)}\,c^k_{m-1}.
\]
If $\d= 1+\frac{2k-l-1}{n+1}$ $($a critical value$)$, there exists a projectively equivariant quantization on $\cS^\d_k$ if and only if~$\l=\frac{1-h}{n+1}$ with $k-l< h \leq k$ $($resonance$)$, and it is given by $\cN\circ\left(\sum_{m=0}^{k}c^k_m D^m\right)$, the coefficients $c^k_m$ being defined as above, except $c^k_l$ which is free.
\end{Theorem}

\subsection{Critical values and cohomology}\label{section3.2}
In light of Theorem \ref{Thmgequi}, the projectively equivariant quantization exists and is unique on~$\cS^\d_k$ if and only if there is no projectively invariant dif\/ferential operator acting on this space.
From preceding considerations, the only candidates are powers of the divergence $D$, and  equation~\eqref{Dproj} shows that $D^l$ is projectively invariant on~$\cS^\d_k$ if and only if $\d=1+\frac{2k-l-1}{n+1}$. Thus, we recover exactly the statement of Theorem \ref{QEP}, and are able to interpret the critical values of $\d$ in terms of existence of projectively invariant operators on~$\cS^\d_k$. Following Lecomte \cite{Lec00}, this can be stated in cohomological terms. Let us introduce the $1$-cocycle
\[
\g(X)=\frac{1}{n+1}\partial_i(\Div X)\partial_{p_i}D^{l-1},
\]
whose class lies in $H^1(\Sl_{n+1},\Hom(\cS^\d_k,\cS^\d_{k-l}))$. It satisf\/ies the equality $\g_l=\left(\ell_{k-l}(\l)c^k_{k-l+1}\right)\g$.  Since $\g$ vanishes on $\aff(n,\bbR)$ and $\Hom_\aff(\cS^\d_k,\cS^\d_{k-l})$ is generated by $D^l$, we deduce from \eqref{Dproj} that $\g$ def\/ines a nontrivial $1$-cocycle if and only if $D^l$ is projectively invariant, i.e.\ $\d= 1+\frac{2k-l-1}{n+1}$.
Consequently, for this critical value of~$\d$, we get an obstruction to the existence of a projectively equivariant quantization except when $\g_l$ vanishes. This occurs when
\begin{gather}\label{LDLd0}
\forall\, X\in\Sl_{n+1},\qquad \big(\LDX-\LdX\big)|_{\cS^\d_{h}}=0 ,
\end{gather}
or equivalently when $\ell_{h-1}(\l)=0$ for some $k-l< h \leq k$, giving the resonant values $(\l,\l+\d)$.

\subsection{Resonances and projectively invariant operators}\label{section3.3}
Here we provide an interpretation of the resonant values of $\l$ in terms of projectively invariant operators acting on the space of densities $\cF^\l$.
\begin{Theorem}\label{ThmProj}
The projectively equivariant quantization $\cQlm:\cS^\d_k\otimes\cF^\l\rightarrow\cF^\m$ exists and is unique except when there is a non-zero projectively invariant differential operator: $\cS^\d_k\rightarrow\cS^\d_{k-l}$. In that case, $\cQlm$ exists if and only if there is a projectively invariant differential operator of order~$h$ from $\cF^\l$ to sections of a homogeneous bundle, with $k-l<h\leq k$.
\end{Theorem}
\begin{proof}
The f\/irst statement has just been proven. The second one relies on the characterization by equation \eqref{LDLd0} of the resonant values $(\l,\l+\d)$. Indeed, this relation provides a way to generate projectively invariant operators.
\begin{Lemma}\label{lemB}
Let $(\LDX-\LdX)|_{\cS^\d_k}=0$ for all $X\in\fkg$, a Lie subalgebra of $\Vect(\bbR^n)$, and let $\mathcal{B}$ be the space of sections of a homogeneous bundle. If there exists a $\fkg$-invariant element in $\cS^\d_k\otimes\mathcal{B}$, then its image by normal ordering is a $\fkg$-invariant differential operator: $\cF^\l\rightarrow\cF^\m\otimes\mathcal{B}$ of order~$k$.
\end{Lemma}
\begin{proof}
The expression $\LDX-\LdX$ is the same if, on the one hand, $\Ld$ is the action on~$\cS^\d_k$ and~$\LD$ is the one on $\D^{\l,\m}_k$, and if, on the other hand, $\Ld$ is the action on $\cS^\d_k\otimes\mathcal{B}$ and $\LD$ is the action on the dif\/ferential operators from $\cF^\l$ to $\cF^\m\otimes\mathcal{B}$ of order $k$.\end{proof}

Consequently, we want to obtain $\Sl_{n+1}$-invariant elements. Using Weyl's theory of invariants of $\gl_n(\bbR)$ ($\subset\Sl_{n+1}$), we have no choice but to take $\mathcal{B}=\bcS^{-\d}=(\bcS(\bbR^n),\bLd)$, the module of functions on $TM$ polynomial in the f\/iber variables tensored with $(-\d)$-densities. Denoting the f\/iber coordinates by $\bar{p}^i$, the $\gl_n(\bbR)$-invariant elements are ${(\bar{p}^ip_i)^h\in\cS_h^\d\otimes\bcS_h^{-\d}}$ for $k\in\bbN$. They are obviously $\Sl_{n+1}$-invariant, and by Lemma~\ref{lemB} and equation~\eqref{ProjLDLd}, they gives rise for $\ell_{h-1}(\l)=0$ to a projectively invariant dif\/ferential operator
\begin{gather}\label{barGproj}
\bar{G}^{h}: \ \cF^\l\rightarrow \bcS^\l_{h},
\end{gather}
with $\bar{G}=\bar{p}^i\partial_i$. As shown by straightforward computations, this is the only projectively inva\-riant dif\/ferential operator with principal symbol $(\bar{p}^ip_i)^{h}$. Thus, there is a projectively invariant dif\/ferential operator of order $h$ on $\cF^\l$ if and only if $\l=\frac{1-h}{n+1}$, and the theorem is proved.
\end{proof}

Now we make concrete the existence of projectively equivariant quantization for resonances, by constructing it from the projectively invariant operators \eqref{barGproj}. Let us denote by $\D(\bcS^\l_k,\cF^\m)$ the space of dif\/ferential operators from $\bcS^\l_k$ to $\cF^\m$, which is isomorphic to $\Dlm\otimes\cS_k^0$ as $\Vect(\bbR^n)$-module. The $\Vect(\bbR^n)$-invariant $1$-form $\a=dx^i\partial_{\bar{p}^i}$ acts by interior product on $\Vect(\bbR^n)$ and vanishes on $\cS^\d_k$. Thus, it extends as a derivation, denoted by $\iota_\alpha$, on the algebra underlying $\Dlm\otimes\cS^0$ and gives rise for any integer $j\leq k$ to the following morphism of $\Vect(\bbR^n)$-modules,
\begin{gather}\label{NbarN}
 (\iota_\alpha)^j: \ \D^{\l,\m}_{k}\rightarrow\D_{k-j}\big(\bcS^\l_{j},\cF^\m\big),
\end{gather}
with kernel $\D^{\l,\m}_{j-1}$. From projective invariance of the operators given in \eqref{barGproj}, we deduce the following proposition.
\begin{Proposition}\label{Propprojinv}
Let $\d=1+\frac{2k-l-1}{n+1}$ and $\l=\frac{1-h}{n+1}$, where $k-l< h\leq k$. Then, the partial projectively equivariant quantization $\cQ^{\l,\m}_{k,h}=\cN\circ(\Id+\dots+\phi_{h-1})$ exists and induces
the following commutative diagram of $\Sl_{n+1}$-modules,
\[
\xymatrix{
\cS^\d_k\otimes\cF^\l\ar[rrrrdd]_{} \ar[rr]^{\cQ^{\l,\m}_{k,h}\otimes \Id\qquad}& & \D^{\l,\m}_k/\D^{\l,\m}_{h-1}\otimes\cF^\l \ar[rr]^{(\iota_\alpha)^{h}\otimes \bar{G}^{h}\quad} && \D_{h-1}(\bcS^\l_{h},\cF^\m)\otimes \bcS^\l_{h}\ar[dd]\\
\\
&&& &\cF^\m
}
\]
This defines a projectively equivariant quantization on $\cS^\d_k$.
\end{Proposition}

\section{Conformally equivariant quantization and invariant operators}\label{section4}
The aim of this section is to obtain a full characterization of existence of the conformally equivariant quantization in the spirit of Theorem \ref{ThmProj}. The task is harder in this case since the spaces $D_{\fkg_{-1}\oplus\fkg_0}(\cS^\d_{k,s},\cS^\d_{l,t})$ are generically multi-dimensional (see \eqref{DecompoS} for notation $\cS^\d_{*,*}$).

\subsection{The conformal Lie algebra}\label{section4.1}
Let $(x^i)$ denote the cartesian coordinates on $\bbR^n$ and $\eta$ be the canonical f\/lat metric of signature~$(p,q)$. The Lie algebra of conformal Killing vector f\/ields on $(\bbR^n,\eta)$,  denoted by $\cf$ for short, is isomorphic to $\ro(p+1,q+1)$. As a subalgebra of polynomial vector f\/ields it admits a~gradation $\cf=\cf_{-1}\oplus\cf_{0}\oplus\cf_{1}$ by the degree of the vector f\/ield components, this is a particular case of IFFT-algebra. The Lie subalgebra $\cf_{-1}$ consists of the translations and $\cf_0$ consists of the linear conformal transformations. Thus, $\cf$ is generated by their sum $\ce=\cf_{-1}\oplus\cf_0$, and by the following conformal inversions in $\cf_1$,
\begin{gather}\label{Xi}
X_i=x_jx^j\partial_i-2x_ix^j\partial_j,
\end{gather}
where $x_i=\eta_{ij}x^j$ and $i=1,\ldots,n$. As $\ce\subset\aff(n,\bbR)$, we have $\LDX=\LdX$ for every $X\in\ce$, and the lack of conformal equivariance for the normal ordering is described by, see \cite{DLO99},
\begin{gather}\label{LDLd}
\LDi-\Ldi=-p_iT+2(\cE+n\l)\partial_{p^i},
\end{gather}
where $X_i$ is the conformal inversion \eqref{Xi} and $T=\eta_{ij}\partial_{p_i}\partial_{p_j}$, $\cE=p_i\partial_{p_i}$.
\subsection{Similarity invariant dif\/ferential operators}\label{section4.2}
The f\/irst step is to describe the space of $\ce$-invariant operators  $\D_\ce(\cS^\d_k,\cS^\d_{k'})$ for arbitrary $k$ and~$k'$. Weyl's theory of invariants \cite{Wey97} insures that the algebra of isometry-invariant dif\/ferential operators acting on $\cS^\d$ is generated by
\[
R=\eta^{ij}p_ip_j, \qquad \cE=p_i\partial_{p_i}, \qquad T=\eta_{ij}\partial_{p_i}\partial_{p_j},
\]
corresponding respectively to the metric (or kinetic energy), the Euler operator, the trace operator, and by
\[
G=\eta^{ij}p_i\partial_j, \qquad D=\partial_i\partial_{p_i}, \qquad L=\eta^{ij}\partial_i\partial_j,
\]
corresponding respectively to the gradient, the divergence and the Laplacian. They all become $\ce$-invariant operators if we consider them as operators from $\cS^\d$ to $\cS^{\d'}$ with a certain shift of weight $\d'-\d$ given in the following table,
\begin{gather}\label{tableInv}
\begin{array}{|c|c|c|c|c|c|}
				\hline	
\text{values of}\;n(\d'-\d)	& -2  & 0 & 2\\[2pt]\hline
\ce\text{-invariant operators}	& T  & \cE, \ D &  R, \ G, \ L \\ \hline
\end{array}
\end{gather}
The computation of the action of $X_i\in\cf_1$ on those operators proves that $\cE$, $R$ and $T$ are conformally invariant, and the table \eqref{tableInv} ensures that $\cE$ and $RT$ preserve the shift. The joint eigenspaces of these two operators def\/ine a decomposition of $\cS^\d$ into $\cf$-submodules, which corresponds to the spherical harmonic decomposition in the $p$ variables,
\begin{gather}\label{DecompoS}
\cS^\d=\bigoplus_{k,s\in\bbN,\, 2s\leq k}\cS^\d_{k,s}.
\end{gather}
More precisely, $\cS^\d_{k,s}$ is the space of homogeneous symbols of degree $k$ of the form $P=R^sQ$ with $TQ=0$. 
Let $2s'\leq k'$ be two integers.
Then, each $\ce$- or $\cf$-invariant operator $\cS^\d_{k,s}\rightarrow\cS^\d_{k',s'}$ gives rise to the following commutative diagram of $\ce$- or $\cf$-modules,
\begin{gather}\label{kslt}
\begin{split}
\xymatrix{
\cS^\d_{k,s}\ar[r]\ar[d]_{T^s} & \cS^\d_{k',s'}\\
\cS^{\d-\frac{2s}{n}}_{k-2s,0}\ar[r] & \cS^{\d-\frac{2s'}{n}}_{k'-2s',0}\ar[u]_{R^{s'}}
}
\end{split}
\end{gather}
We write $G_0$, $D_0$, $L_0$ for the restriction and corestriction of the operators $G$, $D$, $L$ to $\ker T$.
Using the relation $[D,G]=L$, we obtain from~\eqref{tableInv} and~\eqref{kslt} that $\D_\ce(\cS^\d_{k,s},\cS^\d_{k',s'})$ is linearly generated by the monomials $R^{s'}G_0^gL_0^{\ell}D_0^dT^{s}$, such that $g+\ell=s-s'$ and $d+g+2\ell=k-k'$. An explicit description for possibly dif\/ferent weights and $s'=s=0$ follows.
\begin{Proposition}\label{Homce}
Let $k$, $k'$ be integers, $\d,\d'\in\bbR$, and define $\ell=\frac{n}{2}(\d'-\d)-\max(k'-k,0)$. The space $\D_\ce(\cS^\d_{k,0},\cS^{\d'}_{k',0})$ is nontrivial only if $\ell$ is a non-negative integer, and then
\[
\D_\ce\big(\cS^\d_{k,0},\cS^{\d'}_{k',0}\big)=
\begin{cases}
(G_0)^{k'-k}\big\langle  L_0^\ell,G_0L_0^{\ell-1}D_0,\ldots,G_0^\ell D_0^\ell\big\rangle\quad & \text{for } k'-k\geq 0, \\
\big\langle  L_0^\ell,G_0L_0^{\ell-1}D_0,\ldots,G_0^\ell D_0^\ell\big\rangle(D_0)^{k-k'}\quad & \text{otherwise.}
\end{cases}
\]
\end{Proposition}
From this we get a diagrammatic representation of the spaces $\D_\ce(\cS^\d_{k,0},\cS^{\d'}_{k',0})$, one path between two spaces corresponding to a one-dimensional subspace of $\ce$-invariant operators,
\begin{gather}\label{DiagDG}
\begin{split}
\xymatrix{
& & \cS^\d_{k,0}\ar[dl]_{G_0}\ar[dr]^{D_0} && \\
& \cS^{\d+\frac{2}{n}}_{k+1,0}\ar[dl]_{G_0}\ar[dr]^{D_0} && \cS^\d_{k-1,0}\ar[dl]_{G_0}\ar[dr]^{D_0} & \\
&& \cS^{\d+\frac{2}{n}}_{k,0} &&
}
\end{split}
\end{gather}
Since $[D_0,G_0]=L_0$, two paths give rise to independent subspaces if they reach dif\/ferent furthest right column.

\subsection{Classif\/ication of conformally invariant dif\/ferential operators\\ acting on symbols}\label{section4.3}
Resorting to the decomposition \eqref{DecompoS} of $\cS^\d$ into $\cf$-submodules, the classif\/ication of conformally invariant dif\/ferential operators on
$\cS^\d$ amounts to the determination of the $\cf$-invariant dif\/ferential operators between submodules $\cS^\d_{k,s}$ and $\cS^\d_{k',s'}$. The commutative diagram \eqref{kslt} of $\cf$-modules further reduces the quest to the space $\D_\cf(\cS^\d_{k,0},\cS^{\d'}_{k',0})$. The two involved submodules of symbols are modules of sections of irreducible homogeneous bundles associated to the principal f\/iber bundle over $(\bbR^n,\eta)$ of conformal linear frames. Consequently, conformally invariant operators correspond to morphisms of generalized Verma modules. Their classif\/ication has been performed by Boe and Collingwood~\cite{BCo85a,BCo85b}, see also~\cite{ESl97} for a clear summary. This has been translated as a~classif\/ication of conformally invariant dif\/ferential operators by Eastwood and Rice for \mbox{$n=4$}~\cite{ERi87}. Instead applying such an heavy machinery to our peculiar case, we prefer to proceed in a~direct and elementary way.

We compute the space $\D_\cf(\cS^\d_{k,0},\cS^{\d'}_{k',0})$ of conformally invariant operators by searching the elements of $\D_\ce(\cS^\d_{k,0},\cS^{\d'}_{k',0})$ which are invariant under the action of a certain $X\in\cf_1$.
Concerning the $\ce$-invariant operators $D_0^d$, $G_0^g$ and $L_0^\ell$, with source space $\cS^\d_{k,0}$, the action of the conformal inversion $X_i$ follows from \cite{DLO99},
\begin{gather}\label{D0Xi}
\big[D_0^d,\Ldi\big]  = 2d(2k-d-1+n(1-\d))\partial_{p^i}D_0^{d-1},\\  \label{G0Xi}
\big[G_0^g,L_{X_i}^*\big]  = -2g(n\d+g-1)\pi_0p_iG_0^{g-1},\\  \label{L0Xi}
\big[L_0^\ell,L_{X_i}^*\big]  = 2\ell \big[(2(k-\ell)+n(1-2\d))\partial_i+2(G_0\partial_{p^i}-\pi_0p_iD_0)\big]L_0^{\ell-1},
\end{gather}
where $L_{X_i}^*$ is def\/ined in \eqref{Homfkh} and $\pi_0$ is the conformally invariant projection $\pi_0:\cS^\d\rightarrow\ker T$. Since the latter operators generate the spaces $\D_\ce(\cS^\d_{k,0},\cS^{\d'}_{k',0})$ with arbitrary $k'$, $\d'$, we deduce~that,
\begin{gather}\label{Egdl}
[\D_\ce(\cS^\d_{k,0},\cS^{\d'}_{k',0}),L_{X_i}^*]\subset E_D\partial_{p^i}\oplus \pi_0 p_iE_G\oplus\partial_iE_L,
\end{gather}
with $E_G$, $E_D$, $E_L$ the maximal vector spaces generated by the three operators $G_0$, $D_0$, $L_0$, and such that
$G_0E_G, E_DD_0, L_0E_L\subset\D_\ce(\cS^\d_{k,0},\cS^{\d'}_{k',0})$. The independence of the monomials $G_0^gD_0^dL_0^{\ell}$
for dif\/ferent exponents $\ell$ and the commutative diagram \eqref{kslt} lead to the next two results.
\begin{Lemma}\label{lemConfInv}
Let $k$, $k'$ be non-negative integers, $\d,\d'\in\bbR$, and let us define $j=\frac{n}{2}(\d'-\d)$. The space of conformally invariant operators
$\D_\cf(\cS^\d_{k,0},\cS^{\d'}_{k',0})$ is either trivial or of dimension~$1$, generated by
\begin{itemize}\itemsep=0pt
\item $D_0^d$ if $k-k'=d$, $j=0$ and $\d=1+\frac{2k-d-1}{n}$,
\item $G_0^g$ if $k'-k=g$, $j=g$ and $\d=\frac{1-g}{n}$,
\item $\ccL_\ell$ if $k'=k$, $j=\ell$ and $\d=\half+\frac{k-\ell}{n}$,
\end{itemize}
where the operator $\ccL_\ell$ is of the form $L_0^\ell+a_1G_0L_0^{\ell-1}D_0+\dots+a_\ell G_0^\ell D_0^\ell$ for $a_i\in\bbR$. Moreover, the operator $\ccL_\ell$ is not the $\ell^\text{th}$ power of $\ccL_1$. 
\end{Lemma}
\begin{proof}
Let $A\in\D_\cf(\cS^\d_{k,0},\cS^{\d'}_{k',0})$.
Resorting to Proposition \ref{Homce}, up to a constant we have $A=G_0^gBD_0^d$ where $B=L_0^\ell+a_1G_0L_0^{\ell-1}D_0+\dots+a_\ell G_0^\ell D_0^\ell$, for some integers $g$, $d$, $\ell$ and reals $a_1,\ldots,a_\ell$.
As the component in $E_D\partial_{p^i}$ of the higher degree term in $L_0$ of~$[A,L_{X_i}^*]$ vanishes, we get
\[
[D_0^{d},L_{X_i}^*]=0.
\]
From the relation \eqref{D0Xi} we deduce that necessarily $\d=1+\frac{2k-d-1}{n}$ if $d\neq 0$.
As the component in $\pi_0 p_i E_G$ of the higher degree term in $L_0$ of $[A,L_{X_i}^*]$ vanishes, we obtain
\[
[G_0^g,L_{X_i}^*]=0.
\]
From the relation \eqref{G0Xi} we deduce that necessarily $\d+\frac{2\ell}{n}=\frac{1-g}{n}$ if $g\neq 0$.
As a consequence~$B$ must be conformally invariant, and since the component in $\partial_iE_L$ of the higher degree term in~$L_0$ of $[B,L_{X_i}^*]$  vanishes, we are lead to
\[
[L_0^{\ell},L_{X_i}^*]\in \pi_0 p_iE_G\oplus E_D\partial_{p^i}.
\]
From the relation \eqref{L0Xi} we deduce that necessarily $\d=\half+\frac{k-\ell}{n}$ if $\ell\neq 0$. Then, straightforward but lengthy computations show that there exist
unique reals $a_1,\ldots,a_\ell$ such that $L_0^{\ell}+a_1G_0L_0^{\ell-1}D_0+\dots+a_\ell G_0^\ell D_0^\ell$ is conformally invariant, and the expression of $a_1$ ensures that $\ccL_\ell\neq(\ccL_1)^\ell$.
The three values found for $\d$ are incompatible two by two, hence we get that among the exponents~$g$,~$d$,~$\ell$ one at most is non-vanishing. The result follows.
\end{proof}

\begin{Theorem}\label{ConfInv}
Let $k\geq 2s$ and $k'\geq 2s'$ be integers, and $\d,\d'\in\bbR$. The space of conformally invariant differential operators
$\D_\cf(\cS^\d_{k,s},\cS^{\d'}_{k',s'})$ is either trivial or of dimension $1$. In the latter case $j=\frac{n}{2}(\d'-\d)$ is an integer and the space is generated by
\begin{itemize}\itemsep=0pt
\item $R^{s'}D^dT^s$, if $s'-s=j$, $k-k'=d-2j$ and $\d=1+\frac{2(k-s)-d-1}{n}$,
\item $R^{s'}G_0^gT^s$, if $g+s'-s=j$, $k-k'=s-s'-j$ and $\d=\frac{2s+1-g}{n}$,
\item $R^{s'}\ccL_\ell T^s$, if $\ell+s'-s=j$, $k-k'=2(\ell-j)$ and $\d=\half+\frac{k-\ell}{n}$.
\end{itemize}
\end{Theorem}
\begin{Remark}
If $n=2,3$, there is also isometric invariants in $\D(\cS,\cS)$ built from the volume form, but the only conformally invariant one is the algebraic operator $p_1\partial_{p_2}-p_2\partial_{p_1}$ for $n=2$.
\end{Remark}
\begin{Remark}
The conformally invariant operator $G_0^g$ is the $g$-generalized conformal Killing ope\-ra\-tor, whose kernel is the space of $g$-generalized conformal Killing tensors~\cite{NPr90}. The conformally invariant operator $\ccL_\ell$ is the generalization of $\ell^{\text{th}}$ power of the Laplacian to trace free symbols, its curved analog, for $\ell=1$, has been obtained in~\cite{Wun86}.
\end{Remark}

\subsection{Critical values and cohomology}\label{section4.4}
From Theorem \ref{ConfInv}, classifying the conformally invariant dif\/ferential operators on $\cS^\d$, and Theorem~\ref{Thmgequi}, characterizing the critical values of the shift $\d$ in terms of invariant dif\/ferential operators, we recover the following theorem of Silhan.
\begin{Theorem}[\cite{Sil09}]\label{ThmRes}
Let $k\geq 2s$ and $l$ be three integers and $\d\in\bbR$. The conformally equivariant quantization exists and is unique on $\cS^\d_{k,s}$ if and only if there is no conformally invariant differential operator $\cS^\d_{k,s}\rightarrow\cS^\d_{k-l}$. This means $\d\notin I_{k,s}$, where the set of critical values is of the following form: $I_{k,s}=I_{k,s}^D\amalg \big(I_{k,s}^G\cup I_{k,s}^L\big)$ and
\begin{gather} \nonumber
I^D_{k,s} =\left\{1+\frac{2(k-s)-d-1}{n}\left|\; d\in\llbracket 1,k-2s\rrbracket \right.\right\},\\
I^G_{k,s} =\left\{\frac{2s+1-g}{n}\left|\; g\in\llbracket 1,s\rrbracket \right.\right\},\qquad
I^L_{k,s} =\left\{\half+\frac{k-\ell}{n}\left|\; \ell\in\llbracket 1,s\rrbracket \right.\right\}.\label{Igdl}
\end{gather}
\end{Theorem}
We provide now an alternative characterization of the critical values $\d\in I_{k,s}$ in cohomolo\-gi\-cal terms. The Proposition~\ref{Homce} and the commutative diagram \eqref{kslt} allow to decompose the conformally equivariant quantization $\cN\circ\Qlm$ on $\cS^\d_{k,s}$ as
\[
\Qlm=\sum_{g\leq s}\left(\sum_{d\leq k-2s+g}\phi_{d,g}\right),
\]
with $\phi_{d,g}\in\D_{\ce}(\cS^\d_{k,s},\cS^\d_{k-d-g,s-g})$ by $\ce$-invariance of the normal ordering.
Analogously to \eqref{Egdl}, we introduce the spaces $E_ {d,g}=R^{s'}\left(E_D\partial_{p^i}\oplus \pi_0 p_iE_G\oplus\partial_iE_L\right) T^s$,
which are the sum of three maximal subspaces of invariant operators such that
$R^{s'} E_DD_0T^s,\ldots\subset\D_\ce(\cS^{\d}_{k,s},\cS^{\d}_{k-d-g,s-g})$. From Proposition \ref{Homce}, we easily get that $\dim E_{d,0}=\dim E_{0,g}=1$, and denoting $\min (d,g)$ by $\ell$, we also obtain that $\dim E_{d,g}=3\ell+1$ if $g\neq d$ and $\dim E_{\ell,\ell}=3\ell$.

\begin{Lemma}\label{lemEcocycle}
Let $X_i\in\cf_1$. The map $\g\mapsto\g(X_i)$ establishes a linear isomorphism between $E_{d,g}$ and the space of $1$-cocycles vanishing on $\ce$ and with values in $\Hom(\cS^\d_{k,s},\cS^{\d}_{k-d-g,s-g})$.
\end{Lemma}
\begin{proof}
Let $X_i$ and $X_j$ be inversions given by \eqref{Xi}.
From the cocycle relation \eqref{1cocycle} and the vanishing of $\g$ on $\ce$, we get that $\eta^{ij}\g(X_j)\partial_i$ is $\ce$-invariant. This proves that the image of the considered map lies in $E_{d,g}$. Let $A\in E_{d,g}$ and $X_{ji}=x_j\partial_i-x_i\partial_j\in\cf_0$. The latter satisf\/ies $[X_{ji},X_i]=X_j$, and allows to def\/ine a linear map on $\cf$ by: $\g(X)=0$ if $X\in\ce$, $\g(X_i)=A$ and $\g(X_j)=[L^\d_{X_{ji}},A]$ if $j\neq i$. This map is a $1$-cocycle satisfying the required properties.
\end{proof}
\begin{Proposition}\label{Cocyclegdl}
The relative cohomology space $\rH^1(\cf,\ce;\Hom(\cS^\d_{k,s},\cS^{\d}))$ changes of dimension exactly for the critical values $\d\in I_{k,s}$. Then it rises by one, or two if $\d\in I_{k,s}^G\cap I_{k,s}^L$. In particular, the $1$-cocycles
\begin{gather*}
\g_d: \ X \mapsto \partial_i(\Div X)\partial_{p_i}D^{d-1},\qquad
\g_g: \ X \mapsto \eta^{ij} \partial_i(\Div X)R^{s-g}\pi_0p_jG_0^{g-1}T^s,
\end{gather*}
are nontrivial if and only if the operators $D^d$ and $R^{s-g}G_0^{g}T^s$ are $\cf$-invariant.
\end{Proposition}

\begin{proof}
According to the proof of the previous Lemma, if $[\g]\in\rH^1(\cf,\ce;\Hom(\cS^\d_{k,s},\cS^{\d}))$, then $\eta^{ij}\partial_i\g(X_j)$ is $\ce$-invariant. Consequently, the space $\rH^1(\cf,\ce;\Hom(\cS^\d_{k,s},\cS^{\d}))$ is the direct sum of the spaces $\rH^1(\cf,\ce;\Hom(\cS^\d_{k,s},\cS^{\d}_{k-d-g,s-g}))$ for $g\leq s$ and $d\leq k-2s+g$.
The spaces of corresponding $1$-cocycles have just been identif\/ied to the spaces $E_{d,g}$, whose dimension is known and independent of~$\d$. The spaces of corresponding $1$-coboundaries have a dimension equal to that of $\D_\ce(\cS^\d_{k,s},\cS^{\d}_{k-d-g,s-g})$ minus the one of the subspaces of $\cf$-invariant elements. The latter is zero generically, except for the critical values $\d\in I_{k,s}$, hence the result.

For generic $\d$, the two given maps $\g_d$ and $\g_g$ are $1$-coboundaries, proportional to the dif\/feren\-tial of the operators $D^d$ and, respectively, $R^{s-g}G_0^{g}T^s$. Consequently, they def\/ine $1$-cocycles for every $\d\in\bbR$. Since they vanish on $\ce$, they are trivial only if they are of the form $X\mapsto [A,L^\d_X]$, with $A\in\D_\ce(\cS^\d_{k,s},\cS^\d_{k',s'})$ for the adapted values of $k'$ and $s'$.
This space reduces to $0$ if $D^d$ and $R^{s-g}G_0^{g}T^s$ are $\cf$-invariant,
leading to the announced result.
\end{proof}

\subsection{Resonances and conformally invariant dif\/ferential operators}\label{section4.5}

The pairs of weights $(\l,\m)$ for which the conformally equivariant quantization exists are completely known only on symbols of degrees $2$ and $3$ in the momenta \cite{DOv01,Lou01}. Such a pair is called a resonance and in the known cases there are a f\/inite number of resonances for a given critical value of the shift $\d=\m-\l$. We extend those results to the general case and provide a complete classif\/ication for existence of the conformally equivariant quantization.
\begin{Theorem}\label{existQ}
Let $k\geq 2s$ be two integers and $\d\in I_{k,s}$. Restricted to the submodule $\cS^\d_{k,s}$, the existence of the conformally equivariant quantization is equivalent to
\begin{alignat}{5}
& (i)\ && \l=\frac{1-h}{n}, \ &&   k-s-d<h\leq k-s\quad && \text{if} \ \d=1+\frac{2(k-s)-d-1}{n}\in I^D_{k,s}, & \nonumber \\
& (ii)\ &&  \l=\frac{n-2t}{2n}, \  &&  s-g<t\leq s\quad && \text{if} \ \d=\frac{2s+1-g}{n}\in I^G_{k,s}\setminus I^L_{k,s}, &\nonumber \\
& (iii)\ &&  \l=\frac{1-h}{n}, \ &&  k-s-\ell<h\leq k-s \quad && \text{if} \ \d=\half+\frac{k-\ell}{n}\in I^L_{k,s}\setminus I^G_{k,s}, &\nonumber \\
& (iii)'\ && \l=\frac{n-2t}{2n}, \ &&  s-\ell<t\leq s \quad &&\text{if} \ \d=\half+\frac{k-\ell}{n}\in I^L_{k,s}\setminus I^G_{k,s}, &\nonumber \\
& (iv)\ && \l=\frac{n-2t}{2n}, \ && s-\mathrm{min}(\ell,g)<t\leq s \quad && \text{if} \ \d=\half+\frac{k-\ell}{n}=\frac{2g}{n}\in I^L_{k,s}\cap I^G_{k,s}. &\label{lambda}
\end{alignat}
\end{Theorem}
\begin{proof}
Let us sketch the proof.
First, existence of the conformally equivariant quantization translates as the triviality of $1$-cocycles, that we write down in Lemma~\ref{lembbGbbD}. Then, we characterize precisely the conditions for those $1$-cocycles to be trivial in Lemma~\ref{lemconfcoc} and f\/inally prove that, indeed, they are equivalent to \eqref{lambda}.
\begin{Lemma}\label{lembbGbbD}
The conformally equivariant quantization exists on $\cS^\d_{k,s}$ if and only if, for all $g\leq s$ and $d\leq k-2s+g$, there is a map $\phi_{d,g}\in\D_\ce(\cS^\d_{k,s},\cS^\d_{k-d-g,s-g})$ such that
\begin{gather}\label{CocycleTrivconf}
\big[\phi_{d,g},\Ldi\big]=\g_{d,g}(X_i),
\end{gather}
where the right hand side is the corestriction of $(\ccL^{\l,\m}_{X_i}-L^\d_{X_i})\circ(\phi_{d,g-1}+\phi_{d-1,g})$ to the subspace $\cS^\d_{k-d-g,s-g}$. It is
given by
\begin{gather}\label{gammagd}
\g_{d,g}(X_i)=\bbG^{k+1-d-g}_{s+1-g}(\l)\phi_{d,g-1}+\bbD^{k+1-d-g}_{s-g}(\l)\phi_{d-1,g},
\end{gather}
and for any $2s'\leq k'$ the operators $\bbG^{k'}_{s'}(\l):\cS^\d_{k',s'}\rightarrow\cS^\d_{k'-1,s'-1}$ and $\bbD_{s'}^{k'}(\l):\cS^\d_{k',s'}\rightarrow\cS^\d_{k'-1,s'}$ vanish if and only if $\l=\frac{n-2s'}{2n}$ and $\l=\frac{1-(k'-s')}{n}$ respectively.
\end{Lemma}
\begin{proof}
The equation \eqref{CocycleTrivconf} is simply the projection of equation \eqref{phig} on the space $\cS^\d_{k-d-g,s-g}$.
Hence, we just have to get the explicit expression of $\g_{d,g}(X_i)$. For that, we need to decompose the operator $\LDi-\Ldi$ given in \eqref{LDLd}. Let $R^sQ\in\cS^\d_{k,s}$, so that $Q$ is in the kernel of $T$. Using the decomposition $p_iQ=\pi_0(p_iQ)+\frac{2}{\rho_{k-2s+1,1}}R\partial_{p^i}Q$, where $\rho_{k,s}=2s(n+2(k-s-1))$ is the eigenvalue of $RT$ and $\pi_0$ denotes the projection on $\ker T$, we obtain
\begin{gather*}
\big(\LDi-\Ldi\big)(R^sQ) =2s(2n\l+2s-n)R^{s-1}\pi_0(p_iQ)\\
\phantom{\big(\LDi-\Ldi\big)(R^sQ) =}{} +\left(2n\l+2(k-1)+2s\frac{2n\l-2s-n}{n+2(k-2s-1)} \right)R^s\partial_{p^i}Q.
\end{gather*}
Both coef\/f\/icients are af\/f\/ine functions in $\l$ vanishing for the announced values of $\l$, and we have $R^{s-1}\pi_0(p_iQ)\in\cS^\d_{k-1,s-1}$ and $R^s\partial_{p^i}Q\in\cS^\d_{k-1,s}$.
Hence, we get the announced expression~\eqref{gammagd} for $\g_{d,g}(X_i)$.
\end{proof}

\begin{Lemma}\label{lemconfcoc}
Existence of the conformally equivariant quantization on $\cS^\d_{k,s}$ is equivalent to: $\g_{d,0}=0$, $\g_{0,g}=0$ or $\g_{\ell,\ell}(X_i)$ has no component in $\partial_iR^{s-\ell}L_0^{\ell-1}T^s$ if, respectively, $D^d$, $R^{s-g}G_0^gT^s$ or $R^{s-\ell}\ccL_\ell T^s$ is conformally invariant.
\end{Lemma}

\begin{proof}
Using the preceding lemma, we know that the conformally equivariant quantization exists if and only if the $1$-cocycles $\g_{d,g}$ are trivial.
But resorting to the proof of Theorem \ref{Thmgequi}, this is the case except when there is a conformally invariant operator in $\Hom_{\ce}(\cS^\d_{k,s},\cS^\d_{k-d-g,s-g})$. Hence, the previous classif\/ication of such operators shows that the only nontrivial $1$-cocycle is $\g_{d,0}$, $\g_{0,g}$ or $\g_{\ell,\ell}$ if respectively $D^d$, $R^{s-g}G_0^gT^s$ or $R^{s-\ell}\ccL_\ell T^s$ is conformally invariant. Since the previously introduced spaces of $1$-cocycles $E_{d,0}$ and $E_{0,g}$ are unidimensional, Proposition~\ref{Cocyclegdl} implies that $\g_{d,0}$ and $\g_{0,g}$ are trivial if and only if they vanish. It remains to handle the case of~$\g_{\ell,\ell}$.  We suppose that $R^{s-\ell}\ccL_\ell T^s$ is conformally invariant for the shift~$\d$. Let us denote by~$C^{\d'}$ the space of $1$-coboundaries $[\Hom_\ce(\cS^{\d'}_{k,s},\cS^{\d'}_{k-\ell,s-\ell}),L^{\d'}]$. Generically, $\g_{\ell,\ell}$ is a $1$-coboundary, so it pertains to the space $\lim\limits_{\d'\rightarrow\d}C^{\d'}=C^{\d}\oplus\bbR$. The linear form giving the component in $\partial_iR^{s-\ell}L_0^{\ell-1}T^s$ vanishes precisely on the subspace of $1$-coboundaries $C^{\d}$, and thus $\g_{\ell,\ell}$ is trivial if and only if it pertains to its kernel.
\end{proof}

We are ready now to prove the theorem. In each of the four cases in~\eqref{lambda}, i.e.\ for each critical value $\d$, we determine the operators $\bbG^{k'}_{s'}(\l)$ and $\bbD^{k''}_{s''}(\l)$ which should vanish, so that the previous constraints on the $1$-cocycles $\g_{d,0}$, $\g_{0,g}$, $\g_{\ell,\ell}$ are satisf\/ied.

The instances $\d\in I^D_{k,s}$ and $\d\in I_{k,s}^G$ are treated like in the projective case, the concerned space $\Hom_\ce$ being unidimensional.
In the f\/irst case, it exists $d\leq k-2s$ such that $D^d$ is conformally invariant, and Lemma~\ref{lemconfcoc} implies that $\g_{d,0}=0$. Using equation \eqref{gammagd}, we obtain by induction that $\bbD_s^{h+s}(\l)=0$ for some $k-s-d< h\leq k-s$, hence the result (i) of \eqref{lambda}.
If $\d\in I^G_{k,s}$, then it exists $g\leq s$ such that $R^{s-g}G_0^gT^s$ is conformally invariant. Similarly we must have $\g_{0,g}=0$ and then, by induction, $\bbG^{k-s+t}_{t}(\l)=0$ for some $s-g< t\leq s$. The result (ii) of \eqref{lambda} follows.

The instance $\d\in I_{k,s}^L$ is more involved, and relies, as Lemma \ref{lemConfInv}, on independence of mono\-mials $G_0^gD_0^dL_0^{\ell}$ for dif\/ferent exponents $\ell$.
In that case, there exists $\ell\leq s$ such that $R^{s-\ell}\ccL_\ell T^s$ is conformally invariant and Lemma \ref{lemconfcoc} implies that $\g_{\ell,\ell}$ must have no component along $\partial_iR^{s-\ell}L_0^{\ell-1} T^s$. We deduce from the formula \eqref{gammagd} that either $\bbD_{s-\ell}^{k+1-2\ell}=0$ or $\phi_{\ell-1,\ell}$ has no component in $R^{s-\ell}G_0L_0^{\ell-1}T^s$. Then, a straightforward induction proves that all the paths in the following diagram \eqref{DIAG} from $\Id$ to $L_0^\ell$ have a vanishing label on at least one of its arrows. For simplicity we do not write the appropriate powers of $R$ and $T$.
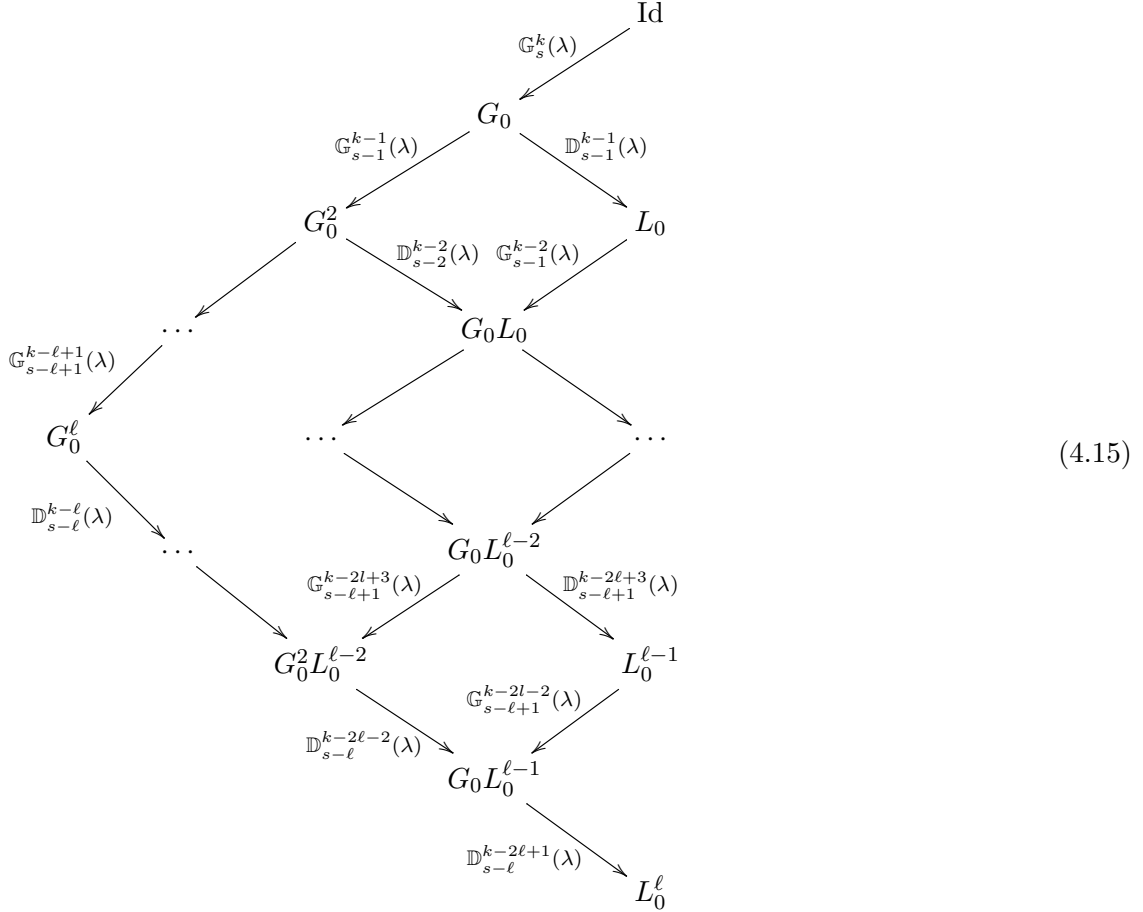
\begin{gather}\label{DIAG}
\begin{split}
\xymatrix{
&&&& \Id \ar[dl]_{\bbG^k_s(\l)}\\
&&& G_0 \ar[dl]_{\bbG^{k-1}_{s-1}(\l)} \ar[dr]^{\bbD^{k-1}_{s-1}(\l)} &\\
&& G_0^2\ar[dl] \ar[dr]^{\bbD^{k-2}_{s-2}(\l)} && L_0 \ar[dl]_{\bbG^{k-2}_{s-1}(\l)}\\
& \cdots \ar[dl]_{\bbG^{k-\ell+1}_{s-\ell+1}(\l)} && G_0L_0 \ar[dl]\ar[dr] &\\
G_0^\ell\ar[dr]_{\bbD^{k-\ell}_{s-\ell}(\l)} && \cdots \ar[dr] && \cdots \ar[dl]\\
& \cdots \ar[dr] && G_0L_0^{\ell-2} \ar[dl]_{\bbG^{k-2l+3}_{s-\ell+1}(\l)} \ar[dr]^{\bbD_{s-\ell+1}^{k-2\ell+3}(\l)}&\\
&& G_0^2L_0^{\ell-2} \ar[dr]_{\bbD^{k-2\ell-2}_{s-\ell}(\l)} && L_0^{\ell-1}\ar[dl]_{\bbG^{k-2l-2}_{s-\ell+1}(\l)}\\
&&& G_0L_0^{\ell-1}\ar[dr]_{\bbD_{s-\ell}^{k-2\ell+1}(\l)} &\\
&&&& L_0^\ell
}
\end{split}
\end{gather}
Since $\bbD^{k-l}_{s-l}(\l)$ and $\bbD^k_s(\l)$ vanish for the same value of $\l$, as well as $\bbG^k_s(\l)$ and $\bbG^{k'}_s(\l)$, we f\/inally get that, necessarily, $\bbG^{k-s+t}_{t}(\l)=0$ or $\bbD^{h+s}_{s}(\l)=0$ for respectively $s-\ell< t\leq s$ and $k-s-\ell<h\leq k-s$. Hence, we end with the results (iii) and (iii)$'$ of \eqref{lambda} for $\d\in I_{k,s}^L\setminus I_{k,s}^G$.

At least, we suppose that $\d\in I^L_{k,s}\cap I^G_{k,s}$. Then, it exists $g,\ell\leq s$ satisfying the relation $\d=\half+\frac{k-\ell}{n}=\frac{2g}{n}$. We deduce that $\l$ must satisfy simultaneously the conditions (ii) and (iii)/(iii)$'$ in~\eqref{lambda} and that leads precisely to~(iv).
\end{proof}

The last part of the proof of Theorem \ref{existQ} shows that each $\l$ given in \eqref{lambda} corresponds to the vanishing of an operator $\bbG$ or $\bbD$, i.e.\ to the vanishing of a certain restriction of $\LD-\Ld$. The Lemma \ref{lemB} allows then to deduce existence of conformally invariant dif\/ferential operators on $\l$-densities, for each resonance $\l$ in \eqref{lambda}. They organize in two families: the powers of Laplacian $\Delta^t:\cF^{\frac{n-2t}{2n}}\rightarrow\cF^{\frac{n+2t}{2n}}$
and the operators $\pi_0\bar{G}^h:\cF^{\frac{1-h}{n}}\rightarrow\bar{\cS}^{\frac{1-h}{n}}_{h,0}$ which are built from the operator $\bar{G}=\bar{p}^i\partial_i$, introduced in \eqref{barGproj}. The latter is dual to the conformal Killing operator $G$ arising in Lemma~\ref{lemConfInv}. By Weyl's theory of invariants or by the general classif\/ication in~\cite{BCo85a,BCo85b}, these are all the conformally invariant dif\/ferential operators on $\l$-densities. Hence, we can reformulate the Theorem~\ref{existQ} along the lines of Kroeske's paradigm. Notice that it is done \`a posteriori, contrary to the Theorem~\ref{Thmgequi}.
\begin{Theorem}
The conformally equivariant quantization $\cQlm:\cS^\d_{k,s}\otimes\cF^\l\rightarrow\cF^\m$ exists and is unique except when there is a conformally invariant differential operator: $\cS^\d_{k,s}\rightarrow\cS^\d_{k-l,s-r}$. Then, $\cQlm$ exists if and only if there is a conformally invariant differential operator from $\cF^\l$ to sections of a homogeneous bundle $\mathcal{V}$, whose principal symbol lies in $\cS^{-\l}_{h,t}\otimes\Ga(\mathcal{V})$ with $\D_{\ce}(\cS_{k,s}^\d,\cS^\d_{h,t})\neq \{0\}$ and $\D_{\ce}(\cS_{k-l,s-r}^\d,\cS^\d_{h,t})= \{0\}$.
\end{Theorem}
Silhan provides in \cite{Sil09} an explicit construction for the conformally equivariant quantization when the shift $\d$ is not one of the critical values \eqref{Igdl}, and also two alternative constructions in the critical cases. In the next two propositions, we prove that they allow to handle all the remaining cases of existence, when $(\l,\l+\d)$ is a resonance \eqref{lambda}.
\begin{Proposition}
In the three cases $(ii)$, $(iii)$, $(iv)$ of \eqref{lambda}, the operator $\Delta^t:\cF^\l\rightarrow\cF^{\l+\frac{2t}{n}}$ is conformally invariant. The commutative diagram
\[
\xymatrix{
\cS^\d_{k,s}\otimes\cF^\l\ar[rrrdd]_{} \ar[rrr]^{a\,T^t\otimes \Delta^t\quad}& && \cS^{\d-\frac{2t}{n}}_{k-2t,s-t}\otimes\cF^{\l+\frac{2t}{n}} \ar[dd]^{\cQ^{\l+\frac{2t}{n},\,\mu}}\\
\\
&& &\cF^\m
}
\]
defines then a conformally equivariant quantization on $\cS^\d_{k,s}$ for a well-chosen constant $a\in\bbR$.
\end{Proposition}
\begin{proof}
The operator $a\,T^t$ above is conformally invariant, and $a$ can be chosen such that the equality $aR^t T^t=\Id$ holds on $\cS^\d_{k,s}$, ensuring the good normalization. Moreover,
the conformally equivariant quantization $\cQ^{\l+\frac{2t}{n},\,\mu}$ is well def\/ined on $\cS^{\d-\frac{2t}{n}}_{k-2t,s-t}$ by Theorem~\ref{ThmRes}.
\end{proof}

Let us introduce a ref\/inement of the f\/iltration of~$\Dlm$, which already appears in \cite{DLO99}.
Namely, we denote by $\D^{\l,\m}_{k,s}$ the subspace of $\D^{\l,\m}_{k}$ given by the image of $\bigoplus_{0\leq s-r\leq k-l}\cS^\d_{l,r}$ under the normal ordering. The expression \eqref{LDLd} of the action of $\cf$ shows that $\D^{\l,\m}_{k,s}$ is in fact a $\cf$-submodule of $\D^{\l,\m}_k$. The subspaces~$\cS^\d_{l,r}$ arising in its def\/inition are characterized equivalently by $\D_\ce(\cS^\d_{k,s},\cS^\d_{l,r})\neq 0$ or, with the notation of Lemma \ref{lembbGbbD}, by their presence in the following tree, analogous to previous ones \eqref{DiagDG} and \eqref{DIAG},
\[
\xymatrix{
&& \cS^\d_{k,s} \ar[dl]_{\bbG^k_s(\l)}\ar[dr]_{\bbD^k_s(\l)}&&\\
& \cS^\d_{k-1,s-1} \ar[dl]_{\bbG^{k-1}_{s-1}(\l)} \ar[dr]_{\bbD^{k-1}_{s-1}(\l)} &&\cS^\d_{k-1,s} \ar[dl]_{\bbG^{k-1}_{s}(\l)} \ar[dr]_{\bbD^{k-1}_{s}(\l)}&\\
\cdots &&\cdots&&\cdots
}
\]
We can restrict the dif\/ferential operators obtained by morphism $(\iota_\alpha)^j$, def\/ined in \eqref{NbarN}, to traceless tensors. The resulting morphism of $\cf$-modules  $(\iota_\alpha)^j_0:\D^{\l,\m}_{k,s}\rightarrow\D^{\l,\m}_{k-h}(\bcS^\l_{h,0},\cF^\m)$ has clearly a nontrivial kernel, which is equal to $\D^{\l,\m}_{j-1+r,r}\supset\D^{\l,\m}_{j-1}$ if $r=\min(s,j-1)$. This leads to the following result, similar to Proposition \ref{Propprojinv} dealing with the projective case.

\begin{Proposition}
In the two cases $(i)$ and $(iii)'$ of \eqref{lambda}, the operator $\pi_0 \bar{G}^h:\cF^\l\rightarrow\bar{S}^\l_{h,0}$ is conformally invariant. The commutative diagram
\[
\xymatrix{
\cS^\d_{k,s}\otimes\cF^\l\ar[rrrrrdd]_{} \ar[rr]^{\cQ^{\l,\m}_{k,h}\otimes \Id\qquad}& & \D^{\l,\m}_{k,s}/\D^{\l,\m}_{k-h,s}\otimes\cF^\l \ar[rrr]^{(\iota_\alpha)^h_0\otimes \pi_0\bar{G}^h} &&& \D^{\l,\m}_{k-h}(\bcS^\l_{h,0},\cF^\m)\otimes \bcS^\l_{h,0}\ar[dd]\\
\\
&&&& &\cF^\m
}
\]
defines then a conformally equivariant quantization on $\cS^\d_{k,s}$.
\end{Proposition}
\begin{proof}
In view of the def\/inition of $h$, there is no obstruction to the existence of the partial conformally equivariant quantization $\cQ^{\l,\m}_{k,h}$. Using the conformal invariance of $\pi_0\bar{G}^h$ and of~$(\iota_\a)^h_0$, we get the result.
\end{proof}

\subsection{Example: conformally equivariant quantization of symbols of degree 3}\label{section4.6}
The conformally equivariant quantization as been explicitly determined up to order $3$ in the momenta by Loubon~Djounga in \cite{Lou01}. We give here for each submodule $\cS^\d_{k,s}$, with $k\leq 3$, the conformally invariant operators acting on it, the corresponding critical values of the shift~$\d$, the resonances $(\l,\l+\d)$ and the associated conformally invariant operators on $\cF^\l$.
$$
\setlength{\extrarowheight}{4pt}
\begin{array}{|c||c|c|c|c|}
\hline
\text{Space} & \text{Op. Inv. on } \cS^\d_{k,s}  & \d & \l & \text{Op. Inv. on } \cF^\l \\[4pt]

\hline
\cS^\d_1 & D & 1 & 0 & \bar{G} \\[4pt]
\hline
\cS^\d_{2,0} & D & \frac{n+2}{n} & -\frac{1}{n} & \bar{G}^2 \\[4pt]
 						& D^2 & \frac{n+1}{n} & -\frac{1}{n},0 & \bar{G}^2,\bar{G} \\[4pt]
\hline
\cS^\d_{2,1} & G_0T & \frac{2}{n} & \frac{n-2}{2n} & \Delta \\[4pt]
 						& L_0T & \frac{n+2}{2n} & \frac{n-2}{2n},0 & \Delta,\bar{G} \\[4pt]
\hline
\cS^\d_{3,0} & D & \frac{n+4}{n} & -\frac{2}{n} & \bar{G}^3 \\[4pt]
 						& D^2 & \frac{n+3}{n} & -\frac{2}{n},-\frac{1}{n} & \bar{G}^3,\bar{G}^2 \\[4pt]
 						& D^3 & \frac{n+2}{n} & -\frac{2}{n},-\frac{1}{n},0 & \bar{G}^3,\bar{G}^2,\bar{G} \\[4pt]
\hline
\cS^\d_{3,1} & RDT & \frac{n+2}{n} & -\frac{1}{n} & \bar{G}^2 \\[4pt]
 						& G_0T & \frac{2}{n} & \frac{n-2}{2n} & \Delta \\[4pt]
 						& \ccL_1T & \frac{n+2}{n} & \frac{n-2}{2n},-\frac{1}{n},0 & \Delta,\bar{G}^2,\bar{G} \\[4pt]
\hline 						
\end{array}
$$
Let $\d$ be one of the above critical values. If we consider the conformally equivariant quantization on the whole space $\cS^\d_3\oplus\cS^\d_2\oplus\cS^\d_1\oplus\cS^\d_0$, as in \cite{Lou01}, then it will exist only for the values of $\l$ appearing in each row where $\d$ is present. Thus, we recover precisely the result of \cite{Lou01}.

\section{Prospects}\label{section5}
We have proven in the setting of IFFT-equivariant quantization, introduced in~\cite{BMa06}, that unique existence of the quantization map is lost if and only if there is an invariant dif\/ferential operator on the space of symbols. The work of Cap and Silhan~\cite{CSi09} allows to trivially extend this result to IFFT-equivariant quantization with values in dif\/ferential operators acting on homogeneous irreducible bundles. It remains then to characterize when equivariant quantization exists but is not unique. We have done it only for projectively and conformally equivariant quantizations, with values in $\Dlm$, in terms of existence of an invariant dif\/ferential operator on the module~$\cF^\l$ of $\l$-densities. In the light of these both examples we propose the following method to link resonances with invariant operators on the source space (e.g.~$\cF^\l$).
\begin{itemize}\itemsep=0pt
\item To each invariant dif\/ferential operator acting on the space of symbols corresponds a nontrivial $1$-cocycle~$\g$.
\item The latter obstructs the existence of the equivariant quantization except for a f\/inite number of values of~$\l$.
\item  For those values, a certain restriction of the operator $\LD-\Ld$ vanishes.
\item An invariant dif\/ferential operator on the source space can then be built, and we obtain all of them in this way.
\item An equivariant quantization can be design from those invariant operators.
\end{itemize}
 We hope that this will lead to the proof of the point $(2)$ in the following conjecture, the point~$(1)$ being Theorem \ref{ThmgequiAHS}.
\begin{Conjecture}
Let $\fkg$ be an IFFT-algebra, $\mathcal{V}$, $\mathcal{W}$ be irreducible homogeneous bundles and  $\mathcal{B}$ a~submodule of $\cS^\d_k(\mathcal{V},\mathcal{W})$. The $\fkg$-equivariant quantization
$
\mathcal{B}\otimes(\Gamma(\mathcal{V})\otimes\cF^\l)\longrightarrow \Gamma(\mathcal{W})\otimes\cF^\m
$
 satisfies the following properties:
\begin{enumerate}\itemsep=0pt
\item[$1)$] it does not exist or is not unique if and only if, for $1\leq l\leq k$, there is a $\fkg$-invariant differential operator: $\mathcal{B}\rightarrow\mathcal{C}$, with $\mathcal{C}$ a submodule of $\cS^\d_{k-l}(\mathcal{V},\mathcal{W})$,
\item[$2)$] for such $\d$, the $\fkg$-equivariant quantization exists if and only if there is a $\fkg$-invariant differential operator on $\Gamma(\mathcal{V})\otimes\cF^\l$, whose principal symbol lies in $\mathcal{A}$ and satisfies the both conditions $\D_{\fkg_{-1}\oplus\fkg_0}(\mathcal{B},\mathcal{A})\neq 0$ and $\D_{\fkg_{-1}\oplus\fkg_0}(\mathcal{C},\mathcal{A})= 0$.
\end{enumerate}
\end{Conjecture}

The next step will be to generalize such a result for $\fkg$-invariant bidif\/ferential pairings, which are bidif\/ferential operators $\Gamma(\mathcal{V})\otimes \Gamma(\mathcal{W})\rightarrow \Gamma(\mathcal{T})$, where $\mathcal{V}$, $\mathcal{W}$, $\mathcal{T}$ are irreducible homogeneous bundles.
Let us discuss the peculiar case of the generalized transvectants or Rankin--Cohen brackets. It has been proved in \cite{ORe03} that there exists a unique conformally invariant bidif\/ferential operator of order $2k$ acting on $\l$- and $\m$-densities, $B_{2k}^{\l,\m}:\cF^\l\otimes\cF^\m\rightarrow\cF^{\l+\m+\frac{2k}{n}}$, if and only if the weights $\l$ and $\m$ do not pertain to the set of exceptional values $\{\frac{n-2k}{2n},\ldots,\frac{n-2}{2n}\}\cup\{\frac{2-2k}{n},\ldots,0\}$.
Remarkably they precisely coincide with the resonances of the conformally equivariant quantization. For those exceptional values, we can construct new  generalized transvectants from the conformally invariant operators $\Delta^\ell$ and $G^g_0$. Explicitly, they are given by $B_{2(k-\ell)}^{\l+\frac{2\ell}{n},\m}\circ (\Delta^\ell\otimes\Id)$ if $\ell\leq k$ and $\l\in\{\frac{n-2k}{2n},\ldots,\frac{n-2}{2n}\}$, and by $\cQlm\circ(G^g_0\otimes\Id)$ if $g\leq 2k-1$, $\mu$ is generic and $\l\in\{\frac{2-2k}{n},\ldots,0\}$. A tight link between bidif\/ferential and dif\/ferential invariant operators shows up one more time. The same kind of idea is used in \cite{BCl12}, where new conformally invariant tri\-linear forms on tensor densities are built from the invariant operators $\Delta^\ell$.
To conclude, Kroeske's paradigm def\/initely asks for deeper investigations.

 \subsection*{Acknowledgements}
 It is a pleasure to acknowledge Christian Duval, Pierre Mathonet and Valentin Ovsienko for fruitful discussions and the referees for suggesting numerous improvements. I thank the Luxem\-bourgian NRF for support via the AFR grant PDR-09-063.

\pdfbookmark[1]{References}{ref}
\LastPageEnding

\end{document}